\documentclass[english]{article}
\usepackage[colorlinks=true,citecolor=black,linkcolor=black,urlcolor=blue]{hyperref}

\newcommand{\arxiv}[1]{\href{http://arxiv.org/abs/#1}{\texttt{arXiv:#1}}}

\def\F{\mathcal{F}}
\def\tr{\textrm}

\def\ex{\tr{ex}}
\def\mc{\mathcal}

\def\prf{\emph{Proof: }}
\def\qed{$\hfill\Box$}

\usepackage{amsmath,amssymb}
\newtheorem{thm}{Theorem}
\newtheorem{lemma}[thm]{Lemma}
\newtheorem{cor}[thm]{Corollary}
\newtheorem{prop}[thm]{Proposition}

\begin{document}

\title{An exact Tur\'an result for tripartite 3-graphs}

\author{Adam Sanitt\thanks{Department of Mathematics, University College London, WC1E 6BT, UK. Email: adam@sanitt.com}\and John
Talbot\thanks{Department of Mathematics, University College London, WC1E 6BT, UK. Email: j.talbot@ucl.ac.uk.}}
\maketitle
\begin{abstract}
Mantel's theorem says that among all triangle-free graphs of a given order the balanced complete bipartite graph is the unique graph of maximum size. We prove an analogue of this result for 3-graphs.   Let $K_4^-=\{123,124,134\}$, $F_6=\{123,124,345,156\}$ and $\mathcal{F}=\{K_4^-,F_6\}$: for $n\neq 5$ the unique $\mathcal{F}$-free 3-graph of order $n$ and maximum size is the balanced complete tripartite 3-graph $S_3(n)$ (for $n=5$ it is $C_5^{(3)}=\{123,234,345,145,125\}$). This extends an old result of Bollob\'as that $S_3(n) $ is the unique 3-graph of maximum size with no copy of $K_4^-=\{123,124,134\}$ or $F_5=\{123,124,345\}$. 
\end{abstract}
\section{Introduction}
If $r\geq 2$ then an \emph{$r$-graph} $G$ is a pair $G=(V(G),E(G))$, where $E(G)$ is a collection of $r$-sets from $V(G)$. The elements of $V(G)$ are called \emph{vertices} and the $r$-sets in $E(G)$ are called \emph{edges}. The number of vertices is the \emph{order} of $G$, while the number of edges, denoted by $e(G)$, is the \emph{size} of $G$.

Given a family of $r$-graphs $\F$, an $r$-graph $G$ is \emph{$\F$-free} if it does not contain a subgraph isomorphic to any member of $\F$. For an integer $n\geq r$ we define the \emph{Tur\'an number of $\F$} to be 
\[
\ex(n,\F)=\max\{e(G):G\tr{ an $\F$-free $r$-graph of order $n$}\}.\]
The related asymptotic \emph{Tur\'an density} is the following limit
(an averaging argument due to Katona, Nemetz and Simonovits \cite{KNS}
shows that it always exists)
\[
\pi\left(\mathcal{F}\right)=\lim_{n\rightarrow\infty}\frac{\textrm{ex}\left(n,\mathcal{F}\right)}{\binom{n}{r}}.
\]

The problem of determining the Tur\'an density is essentially solved
for all 2-graphs by the Erd\H os--Stone--Simonovits Theorem.
\begin{thm}[Erd\H os and Stone \cite{ESt}, Erd\H os and Simonovits \cite{ESi}]
Let $\mathcal{F}$ be a family of 2-graphs. If $t=\min\left\{ \chi(F)\,:\, F\in\mathcal{F}\right\} \geq2$,
then
\[
\pi\left(\mathcal{F}\right)=\frac{t-2}{t-1}.
\]
\end{thm}
It follows that the set of all\emph{ }Tur\'an densities for 2-graphs
is $\{0,1/2,2/3,3/4,\ldots\}$.

There is no analogous result for $r\geq3$ and most progress has been
made through determining the Tur\'an densities of individual graphs
or families of graphs. A central problem, originally posed by Tur\'an,
is to determine $\ex(n,K_{4}^{(3)})$, where $K_{4}^{(3)}=\{123,124,134,234\}$ is
the complete 3-graph of order 4. This is a natural extension of
determining the Tur\'an number of the triangle for 2-graphs, a question
answered by Mantel's theorem \cite{M}. Tur\'an gave a construction that he conjectured to be optimal that has density
$5/9$ but this question remains unanswered despite a great deal of
work. The current best upper bound for $\pi(K_4^{(3)})$ is $0.561666$, given by Razborov \cite{R4}.

A related problem due to Katona is given by considering cancellative hypergraphs. A hypergraph $H$ is \emph{cancellative} if for any distinct edges $a,b\in H$, there is no edge $c\in H$ such that $a\triangle b\subseteq c$ (where $\triangle$ denotes the symmetric
difference). For 2-graphs, this is equivalent to forbidding all triangles.
For a 3-graph, it is equivalent to forbidding the two non-isomorphic
configurations $K_{4}^{-}=\{123,124,134\}$ and $F_{5}=\{123,124,345\}$.

An $r$-graph $G$ is \emph{$k$-partite} if there is a partition of its vertices into $k$ classes so that all edges of $G$ contain at most one vertex from each class. It is \emph{complete $k$-partite} if there is a partition into $k$ classes such that all edges meeting each class at most once are present. If the partition of the vertices of a complete $k$-partite graph is into classes that are as equal as possible in size then we say that $G$ is \emph{balanced}.

Let $S_3(n)$ be the complete balanced tripartite 3-graph of order $n$.
\begin{thm}[Bollob\'as \cite{B}]\label{bol:thm}
For $n\geq 3$, $S_3(n)$ is the unique cancellative 3-graph of order $n$ and maximum size.
\end{thm}
This result was refined by Frankl and F\"uredi \cite{FF} and
Keevash and Mubayi \cite{KM}, who proved that $S_3(n)$ is
the unique $F_5$-free 3-graph of order $n$ and maximum size, for
 $n$ sufficiently large.

The \emph{blow-up} of an $r$-graph $H$ is the $r$-graph $H(t)$ obtained from $H$ by
replacing each vertex $a\in V(H)$ with a set of $t$ vertices $V_{a}$
in $H(t)$ and inserting a complete $r$-partite $r$-graph between any $r$ vertex classes corresponding to an edge in $H$. The following result is an invaluable
tool in determining the Tur\'an density of an $r$-graph that is contained in the blow-ups of other $r$-graphs:
\begin{thm}[Brown and Simonovits \cite{BS}, \cite{BT}]
\label{thm:blowup}If $F$ is a $k$-graph that is contained in a blow-up of every member of a family of $k$-graphs $\mathcal{G}$, then $\pi\left(F\right)=\pi\left(F\cup\mathcal{G}\right)$.
\end{thm}
Since $F_5$ is contained in $K_4^-(2)$, Theorems \ref{bol:thm} and \ref{thm:blowup} imply that $\pi(F_5)=2/9$.

A natural question to ask is which $3$-graphs (that are not subgraphs of blow-ups of $F_5$) also have Tur\'an density $2/9$? Baber and Talbot \cite{BT} considered the 3-graph $F_{6}=\{123,124,345,156\}$,
which is not contained in any blow-up of $F_{5}$. Using Razborov's
flag algebra framework \cite{RF}, they gave a computational
proof that $\pi\left(F_{6}\right)=2/9$. In this paper,
we obtain a new (non-computer) proof of this result. In fact we go further and determine the exact Tur\'an number of $\F=\{F_6,K_4^-\}$. 
\begin{thm}\label{main:thm}If $n\geq3$
then the unique $\mathcal{F}$-free 3-graph with $\textrm{ex}(n,\mathcal{F})$
edges and $n$ vertices is $S_3(n)$ unless $n=5$ in which case
it is $C_{5}^{(3)}$. \end{thm}
As $F_{6}$ is contained in $K_{4}^{-}(2)$, we have the following corollary to Theorem \ref{thm:blowup}.
\begin{cor}\label{thm:f6turan}$\pi\left(F_{6}\right)=2/9$.\end{cor}

\section{Tur\'an number}
\emph{Proof of Theorem \ref{main:thm}:} We use induction on $n$. Note that the result holds trivially for $n=3,4$. For $n=5$ it is straightforward to check that the only $\F$-free 3-graphs with 4 edges are $S_3(5)$, $\{123,124,125,345\}$ and $\{123,234,345,451\}$. Of these the first two are edge maximal while the third can be extended by a single edge to give $C_5^{(3)}$. Thus we may suppose that $n\geq 6$ and the theorem is true for $n-3$.

For $k\geq 2$ let $T_k(n)$ be the $k$-partite Tur\'an graph of order $n$: this is the complete balanced $k$-partite graph. We denote the number of edges in $S_3(n)$ and $T_k(n)$ by $s_3(n)$ and $t_k(n)$ respectively. 
Let $G$ be $\F$-free with $n\geq 6 $ vertices and $\ex(n,\F)$ edges. Since $S_3(n)$ is $\F$-free we have $e(G)\geq s_3(n)$.

The inductive step proceeds as follows: select a special edge $abc\in E(G)$ (precisely how we choose this edge will be explained in Lemma \ref{key:lem} below). For $0\leq i \leq 3$ let $f_i$ be the number of edges in $G$ meeting $abc$ in exactly $i$ vertices. By our inductive hypothesis we have
\begin{equation}\label{ind1:eq}
e(G)=f_0+f_1+f_2+f_3\leq\ex(n-3,\F)+f_1+f_2+1.\end{equation}
Note that unless $n-3=5$ our inductive hypothesis says that $\ex(n-3,\F)=s_3(n-3)$ with equality iff $G-\{a,b,c\}=S_3(n-3)$. For the moment we will assume that $n\neq 8$ and so we have the following bound 
\begin{equation}\label{s3:eq}
e(G)\leq s_3(n-3)+f_1+f_2+1,
\end{equation}
with equality iff $G-\{a,b,c\}=S_3(n-3)$.

Let $V^-=V(G)-\{a,b,c\}$. For each pair $xy\in \{ab,ac,bc\}$ define $\Gamma_{xy}=\{z\in V^-:xyz\in E(G)\}$ and let $\Gamma_{abc}=\Gamma_{ab}\cup \Gamma_{ac}\cup\Gamma_{bc}$ be the \emph{link-neighbourhood} of $abc$. Note that since $G$ is $K_4^-$-free this is a disjoint union, so \[
f_2=|\Gamma_{ab}|+|\Gamma_{ac}|+|\Gamma_{bc}|=|\Gamma_{abc}|.\] 

For  $x\in \{a,b,c\}$ define $L(x)$ to be the \emph{link-graph of $x$}, so $V(L(x))=V^-$ and $E(L(x))=\{yz\subset V^-:xyz\in E(G)\}$. The \emph{link-graph of the edge $abc$} is  the edge labelled graph $L_{abc}$ with vertex set $V^-$ and edge set $L(a)\cup L(b) \cup L(c)$. The label of an edge $yz\in E(L_{abc})$ is $l(yz)=\{x\in \{a,b,c\}:xyz\in E(G)\}$. The \emph{weight} of an edge $yz\in L_{abc}$ is $|l(yz)|$ and the weight of $L_{abc}$ is $w(L_{abc})=\sum_{yz\in L_{abc}}|l(yz)|$. Note that  $f_1=w(L_{abc})$.

By a subgraph of $L_{abc}$ we mean an ordinary subgraph of the underlying graph where the labels of edges are non-empty subsets of the labels of the edges in $L_{abc}$. For example if $xy\in E(L_{abc})$ has $l(xy)=ab$ then in any subgraph of $L_{abc}$ containing the edge $xy$ it must have label $a,b$ or $ab$.

A triangle in $L_{abc}$ is said to be \emph{rainbow} iff all its edges have weight one and are labelled $a,b,c$. Given an edge labelled subgraph $H$ of $L_{abc}$ and an (unlabelled) graph $G$ we say that $H$ is a \emph{rainbow} $G$ if all of the edges in $H$ have weight 1 and all the triangles in $H$ are rainbow.

The following lemma provides our choice of edge $abc$.
\begin{lemma}\label{key:lem}If $G$ is an $\F$-free 3-graph with $n\geq 6$ vertices and $\ex(n,\F)$ edges then there is an edge $abc\in E(G)$ such that
\[
w(L_{abc})+|\Gamma_{abc}|\leq t_3(n-3)+n-3,
\]
with equality iff $L_{abc}$ is a rainbow $T_3(n-3)$ and $\Gamma_{abc}=V^-$.
\end{lemma}

Underlying all our analysis are some simple facts regarding $\mc{F}$-free 3-graphs that are contained in Lemmas \ref{ax:lem} and \ref{struct1:lem}.
\begin{lemma}\label{ax:lem} If $G$ is $\F$-free and $abc\in E(G)$ then the following configurations cannot appear as subgraphs of $L_{abc}$. Moreover any configuration that can be obtained from one described below by applying a permutation to the labels $\{a,b,c\}$ must also be absent.
\begin{itemize}
\item[($F_6$-1)] The triangle $xy,xz,yz$ with $l(xy)=l(xz)=a$ and $l(yz)=b$.
\item[($F_6$-2)] The pair of edges $xy,xz$ with $l(xy)=ab$ and $l(xz)=c$.
\item[($F_6$-3)] A vertex $x\in \Gamma_{ab}$ and edges $xy,yz$ with labels $l(xy)=c$ and $l(yz)=a$.
\item[($F_6$-4)] A vertex $x\in \Gamma_{ab}$ and edges $xy,yz,zw$ with labels $l(xy)=l(zw)=a$ and $l(yz)=b$.
\item[($F_6$-5)] Vertices $x\in \Gamma_{ac},y\in \Gamma_{bc},z\in \Gamma_{ab}$ and the edge $xy$ with label $l(xy)=b$.
\item[($K_4^-$-1)] The triangle $xy,xz,yz$ with $l(xy)=l(xz)=l(yz)=a$.
\item[($K_4^-$-2)] The vertex $x\in \Gamma_{ab}$ and edge $xy$ with label $l(xy)=ab$.
\item[($K_4^-$-3)] The vertices $x,y\in \Gamma_{ab}$ and edge $xy$ with label $l(xy)=a$.
\end{itemize}
\end{lemma}
\begin{lemma}\label{struct1:lem}If $G$ is $\F$-free and  $abc\in E(G)$ then the link-graph and link-neighbourhood satisfy:
\begin{itemize}
\item[(i)] The only triangles in $L_{abc}$ are rainbow.
\item[(ii)] The only $K_4$s in $L_{abc}$ are rainbow.  
\item[(iii)] $L_{abc}$ is $K_5$-free.
\item[(iv)] If $xy\in E(L_{abc})$ has $l(xy)=abc$ then $x$ and $y$ meet no other edges in $L_{abc}$ and $x,y \not\in \Gamma_{abc}$.
\item[(v)] If $V_{abc}^4=\{x\in V^-:\tr{there is a $K_4$ containing $x$}\}$ 
then $\Gamma_{abc}\cap V_{abc}^4=\emptyset$.
\item[(vi)] There are no edges in $L_{abc}$ between $\Gamma_{abc}$ and $V_{abc}^4$.
\item[(vii)] If $x\in V_{abc}^4$ then $|l(xy)|\leq 1$ for all $y\in V^-$.
\item[(viii)] If $x\in \Gamma_{ac}$, $y\in \Gamma_{bc}$ and $l(xy)=ab$, then $\Gamma_{bc}=\emptyset$. Moreover, if $xz\in E(L_{abc})$ with $z\neq y$ then $z\not\in \Gamma_{abc}$ and $l(xz)=a$, while if $yz\in E(L_{abc})$ with $z\neq x$ then $z\not\in \Gamma_{abc}$ and $l(yz)=b$.
\item[(ix)] If $xy,xz\in E(L_{abc})$, $l(xy)=ab$ and $z\in \Gamma_{abc}$ then $|l(xz)|\leq 1$.
\end{itemize}
\end{lemma} 

We also require the following identities, that are easy to verify. 
\begin{lemma}\label{id:lem}
If $n\geq k \geq 3$ then
\begin{itemize}
\item[(i)] $s_3(n)=s_3(n-3)+t_3(n-3)+n-2$.
\item[(ii)]$t_3(n)=t_3(n-3)+2n-3$.
\item[(iii)]$t_3(n)=t_3(n-2)+n-1+\lfloor n/3\rfloor$.
\item[(iv)]$t_k(n)=t_k(n-1)+n-\lceil n/k\rceil.$
\end{itemize}
\end{lemma}
Let $abc\in E(G)$ be a fixed edge given by Lemma \ref{key:lem}.

By assumption $e(G)\geq s_3(n)$ so Lemma \ref{id:lem}  (i) and Lemma \ref{key:lem} together with the bound on $e(G)$ given by (\ref{s3:eq}) imply that $e(G)=s_3(n)$ and hence $G-\{a,b,c\}=S_3(n-3)$, $L_{abc}$ is a rainbow $T_3(n-3)$ and $\Gamma_{abc}=V^-$. 
To complete the proof we need to show that $G=S_3(n)$. First note that since $L_{abc}$ is a rainbow $T_3(n-3)$ and $\Gamma_{abc}=V^-$, Lemma \ref{struct1:lem} (i) and Lemma \ref{ax:lem}($F_6$-3) imply that no vertex in $\Gamma_{ab}$ is in an edge with label $c$ and similarly for $\Gamma_{ac},\Gamma_{bc}$. Hence $L_{abc}$ is the complete tripartite graph with vertex classes $\Gamma_{ab}$, $\Gamma_{ac}$ and $\Gamma_{bc}$ and the edges between any two parts are labelled with the common label of the parts (e.g.~all edges from $\Gamma_{ab}$ to $\Gamma_{ac}$ receive label $a$). So $L_{abc}$ is precisely the link graph of an edge $abc\in S_3(n)$.

In order to deduce that $G=S_3(n)$ we need to show that $G-\{a,b,c\}=S_3(n-3)$ has the same tripartition as $L_{abc}$. This is straightforward: any edge $xyz\in E(G-\{a,b,c\})$ not respecting the tripartition of $L_{abc}$ meets one of the parts at least twice. But if $x,y,z \in \Gamma_{ab}$ then $|\Gamma_{ac}|\geq 2$ so let $u\in \Gamma_{ac}$. Setting $a=1,b=2,x=3,y=4,z=5,u=6$ gives a copy of $F_6$.  If $x,y\in \Gamma_{ab}$ and $z\in \Gamma_{ac}$ then $a=1,x=3,y=4,z=2$ gives a copy of $K_4^-$.

Hence $G=S_3(n)$ and the proof is complete in the case $n\neq 8$.

For $n=8$ we note that if $G-\{a,b,c\}$ is $F_5$-free then Theorem \ref{bol:thm} implies that the result follows as above, so we may assume that $G-\{a,b,c\}$ contains a copy of $F_5$. In this case it is sufficent to show that $e(G)\leq 17 <18 =s_3(8)$.

If $V(G-\{a,b,c\})=\{s,t,u,v,w\}$ then we may suppose that $stu, stv,uvw, abc \in G$. Since $G$ is $K_4^-$-free it does not contain $suv$ or $tuv$. Moreover it contains at most 3 edges from $\{u,v,w\}^{(2)}\times\{a,b,c\}$ and at most 5 edges from $\{s,t,u,v,w\}\times \{a,b,c\}^{(2)}$. Since $G$ is $F_6$-free it contains no edges from $\{s,t\}\times \{w\}\times\{a,b,c\}$. 

The only potential edges we have yet to consider are those in  $\{st,su,tu,sv,tv\}\times\{w,a,b,c\}$. Since $G$ is $K_4^-$-free it contains at most 2 edges from $std,sud,tud,svd,tvd$, for any $d\in \{w,a,b,c\}$. Moreover, since $G$ is $F_6$-free, if it contains 2 such edges for a fixed $d$ then it can contain at most 3 such edges in total for the other choices of $d$. Hence at most 5 such edges are present.

Thus in total $e(G)\leq 4+ 3 +5 +5=17$, as required.
\qed

In order to prove Lemma \ref{key:lem} we first need an edge with large link-neighbourhood.
\begin{lemma}\label{edge:lem}If $G$ is $K_4^-$-free 3-graph of order $n$ with $s_3(n)$ edges, then there is an edge $abc\in E(G)$ with $|\Gamma_{abc}|\geq n-\lfloor n/3\rfloor-3$.
\end{lemma}
\emph{Proof of Lemma \ref{edge:lem}:} 
Let $G$ be $K_4^-$-free with $n$ vertices and $s_3(n)$ edges. For $x,y\in V(G)$ let $d_{xy}=|\{x:xyz\in E(G)\}$. If $uvw\in E(G)$ then $\Gamma_{uvw}=\Gamma_{uv}\cup\Gamma_{uw}\cup \Gamma_{vw}$ is a union of pairwise disjoint sets and $|\Gamma_{uvw}|=d_{uv}+d_{uw}+d_{vw}-3$. Thus if the lemma fails to hold then for every edge $uvw\in E(G)$ we have $d_{uv}+d_{uw}+d_{vw}\leq n-\lfloor n/3\rfloor -1$. Note that since $\sum_{xy\in \binom{V}{2}}d_{xy}=3e(G)$, convexity implies that
\[
e(G)(n-\left\lfloor \frac{n}{3}\right\rfloor -1)\geq\sum_{uvw\in E(G)}d_{uv}+d_{uw}+d_{vw}=\sum_{xy\in \binom{V}{2}}d_{xy}^2\geq \frac{9e^2(G)}{\binom{n}{2}}.\]
Thus
\[
e(G)\leq \frac{1}{18}n(n-1)(n-\lfloor n/3\rfloor -1).\]
But it is easy to check that this is less than $s_3(n)$. \qed

Our next objective is to describe various properties of the link-graph $L_{abc}$ and link-neighbourhood $\Gamma_{abc}$.

Lemma \ref{struct1:lem} (v) allows us to partition the vertices of $L_{abc}$ as $V^-=\Gamma_{abc}\cup V_{abc}^4\cup R_{abc}$, where $V_{abc}^4=\{x\in V^-:\tr{there is a $K_4$ containing $x$}\}$ and $R_{abc}=V^--(\Gamma_{abc}\cup V_{abc}^4)$. To prove Lemma \ref{key:lem} we require the following result to deal with the part of $L_{abc}$ not meeting any copies of $K_4$.
\begin{lemma}\label{work:lem} Let $H$ be a subgraph of $L_{abc}$ with $s\geq 3$ vertices satisfying $V(H)\cap V_{abc}^4=\emptyset$. If $H_\Gamma=V(H)\cap \Gamma_{abc}$ and $|H_\Gamma|\geq s-\lfloor s/3\rfloor -1$ then \[
w(H)+|H_\Gamma|\leq t_3(s)+s,\] with equality iff $H_\Gamma=V(H)$ and $H$ is a rainbow $T_3(s)$.
\end{lemma}
\emph{Proof of Lemma \ref{key:lem}:} Let $G$ be $\F$-free with $n\geq 6$ vertices and $\ex(n,\F)$ edges. By Lemma \ref{edge:lem} we can choose  an edge $abc\in E(G)$ such that $|\Gamma_{abc}|\geq n-\lfloor n/3 \rfloor -3$. Let $V^-=\Gamma_{abc}\cup R_{abc}\cup V_{abc}^4$ be the partition of $V^-$ given by Lemma \ref{struct1:lem} (v). If $s=|V^-|$, $j=|\Gamma_{abc}|$, $k=|R_{abc}|$ and $l=|V_{abc}^4|$ then $n-3=s=j+k+l$ and $j\geq s-\lfloor s/3\rfloor-1\geq j+k-\lfloor(j+k)/3\rfloor-1$ . We can apply Lemma \ref{work:lem} to $H=L_{abc}[\Gamma_{abc}\cup R_{abc}]$, to deduce that \[
w(L_{abc}[\Gamma_{abc}\cup R_{abc}])+|\Gamma_{abc}|\leq t_3(j+k)+j+k,\] with equality iff $R_{abc}=\emptyset$ and $L_{abc}[\Gamma_{abc}]$ is a rainbow $T_3(j+k)$. Now if $L_{abc}$ is $K_4$-free then $V_{abc}^4=\emptyset$ and the proof is complete, so suppose there is a $K_4$ in $L_{abc}$. In this case $4\leq |V_{abc}^4|\leq n-3-|\Gamma_{abc}|\leq \lfloor n/3\rfloor$, so $n\geq 12$.

We now need to consider the edges in $L_{abc}$ meeting $V_{abc}^4$. By Lemma \ref{struct1:lem} (iii) we know that $L_{abc}$ is $K_5$-free, while Lemma \ref{struct1:lem} (vii) says that $V_{abc}^4$ meets no edges of weight 2 or 3, so by Tur\'an's theorem $w(L_{abc}[V_{abc}^4])\leq t_4(l)$. 

Lemma \ref{struct1:lem} (vi) implies that there are no edges from $\Gamma_{abc}$ to $V_{abc}^4$ so the total weight of edges between $\Gamma_{abc}\cup R_{abc}$ and $V_{abc}^4$ is at most $kl$.  Thus
\[
w(L_{abc})+|\Gamma_{abc}|\leq t_3(j+k)+j+k+t_4(l)+kl.\] Finally Lemma \ref{k4:lem} with $s=n-3$ implies that
\[
w(L_{abc})+|\Gamma_{abc}|\leq t_3(n-3)+n-3,\] with equality iff $R_{abc}=V_{abc}^4=\emptyset$ and $L_{abc}$ is a rainbow $T_3(n-3)$ as required. \qed
\begin{lemma}\label{k4:lem}If $j,k,l\geq 0$ are integers satisfying $j+k+l=s\geq 5$ and $j\geq s-\lfloor s/3\rfloor -1$ then
\begin{equation}\label{k4:eq}
t_3(j+k)+t_4(l)+j+k+kl\leq t_3(s)+s,\end{equation}
with equality iff $l=0$.
\end{lemma}
\emph{Proof of Lemma \ref{k4:lem}:}
If $l=0$ then the result clearly holds, so suppose that $l\geq 1$, $j+k+l=s\geq 5$ and $j\geq s-\lfloor s/3 \rfloor -1$. Let $f(j,k,l)$ be the LHS of (\ref{k4:eq}).  We need to check that $\Delta(j,k,l)=f(j,k+1,l-1)-f(j,k,l)>0$. Using Lemma \ref{id:lem} (iv) we have
\begin{eqnarray*}
\Delta(j,k,l)&=&j-\lceil (j+k+1)/3\rceil +\lceil l/4\rceil+1\\
& = & j+\lceil l/4\rceil-\lfloor (j+k)/3\rfloor.
\end{eqnarray*}
So it is sufficient to check that $j+l/4>(j+k)/3$. This follows easily from $j\geq s-\lfloor s/3\rfloor -1$,  $k\leq \lfloor s/3\rfloor +1$, $l\geq 1$ and $s\geq 5$.\qed

\emph{Proof of Lemma \ref{work:lem}:}
We prove this by induction on $s\geq 3$. 
The result holds for $s=3,4$ (see the end of this proof for the tedious details) so suppose that $s\geq 5$ and the result holds for $s-2$.

Let $H$ be a subgraph of $L_{abc}$ with $s\geq 5$ vertices satisfying $V(H)\cap V_{abc}^4=\emptyset$. Let $H_\Gamma=V(H)\cap \Gamma_{abc}$ and suppose that $|H_\Gamma|\geq s-\lfloor s/3\rfloor -1$. 

Note that if $H$ contains no edges of weight 2 or 3 then the result follows directly from Tur\'an's theorem and Lemma \ref{struct1:lem} (i), so we may suppose there are edges of weight 2 or 3. With this assumption it is sufficient  to show that 
\[
w(H)+|H_\Gamma|\leq t_3(s)+s-1.\]
By Lemma \ref{id:lem} (iii) this is equivalent to showing that the following inequality holds:
\begin{equation}\label{ind:eq}
w(H)+|H_\Gamma|\leq t_3(s-2)+2s-2+\lfloor s/3 \rfloor
\end{equation} 
\textbf{Case (i):} There exists an edge of weight 3, $l(xy)=abc$.

Lemma \ref{struct1:lem} (iv) implies that $x,y\not\in H_{\Gamma}$ and $x,y$ meet no other edges in $H$, so we can apply the inductive hypothesis to $H'=H-\{x,y\}$ to obtain
\[
w(H)+|H_\Gamma|\leq w(H')+|H'_\Gamma|+3\leq t_3(s-2)+s-2+3.\]
Hence (\ref{ind:eq}) holds as required. So we may suppose that $H$ contains no edges of weight 3.

\textbf{Case (ii):} The only edges of weight 2 are contained in $H_\Gamma$

Let $xy\in E(H)$ have weight 2, say $l(xy)=ab$.
Now Lemma \ref{ax:lem} $(K_4^-$-$2)$ implies that $x,y\not\in \Gamma_{ab}$, while Lemma \ref{ax:lem} $(K_4^-$-$3)$ implies that $x,y$ cannot both belong to $\Gamma_{ac}$ or $\Gamma_{bc}$ so we may suppose that $x\in \Gamma_{ac}$ and $y\in \Gamma_{bc}$. Lemma \ref{struct1:lem} (viii) implies that $x,y$ have no more neighbours in $H_\Gamma$. If $H_\Gamma=V(H)$ then we can apply the inductive hypothesis to $H'=H-\{x,y\}$ to obtain
\[
w(H)+|H_\Gamma|\leq t_3(s-2)+s-2+2+2,\] in which case (\ref{ind:eq}) holds, so suppose $V(H)\neq H_\Gamma$. 

Let $z\in V(H)-H_\Gamma$ be a neighbour of $x$ in $H$ if one exists otherwise let $z$ be any vertex in $V(H)-H_\Gamma$. By our assumption that all edges of weight 2 are contained in $H_\Gamma$, $z$ meets no edges of weight 2. Moreover, by Lemma \ref{struct1:lem} (viii), all edges containing $x$ (except $xy$) have label $b$, so $x$ is not in any triangles in $H$. Hence $x$ and $z$ have no common neighbours in $H$ and so the total weight of edges meeting $\{x,z\}$ is at most $2+1+s-3$ (if $xz$ is an edge) and at most $2+s-2$ otherwise. Applying our inductive hypothesis to $H'=H-\{x,z\}$ we have
\[
w(H)+|H_\Gamma|\leq t_3(s-2)+s-2+1+s,\] and (\ref{ind:eq}) holds.

\textbf{Case (iii):} There is an edge of weight 2 meeting $V(H)-H_\Gamma$.

So suppose that $xy\in E(H)$, $l(xy)=ab$ and $y\not\in H_\Gamma$. Lemma \ref{struct1:lem} (ix) implies that for any $z\in H_\Gamma$ we have $|l(xz)|$, $|l(yz)|\leq 1$. Let $\gamma_{xy}=|\{x,y\}\cap H_\Gamma|\leq 1$. Thus, since $xy$ is not in any triangles, the total weight of edges meeting $\{x,y\}$ is at most \[
2+s-2+|V(H)-H_\Gamma|-(2-\gamma_{xy}).\]
Applying the inductive hypothesis to $H'=H-\{x,y\}$ we have
\[
w(H)+|H_\Gamma|\leq t_3(s-2)+s-2+s+s-|H_\Gamma|-2+2\gamma_{xy},\] with equality holding only if $|H'_\Gamma|=s-2$.
Now $|H_\Gamma|\geq s-\lfloor s/3\rfloor -1$ implies that
\begin{equation}\label{strict:eq}
w(H)+|H_\Gamma|\leq t_3(s-2)+2s-3+\lfloor s/3\rfloor +2\gamma_{xy},\end{equation}
 with equality only if $|H'_\Gamma|=s-2$ and $|H_\Gamma|=s-\lfloor s/3\rfloor -1$. If $\gamma_{xy}=0$ then (\ref{ind:eq}) holds as required, so suppose $\gamma_{xy}=1$. In this case (\ref{ind:eq}) holds, unless (\ref{strict:eq}) holds with equality. But if (\ref{strict:eq}) is an equality then $|H_\Gamma|=|H'_\Gamma|+1=s-1$, while $|H_\Gamma|=s-\lfloor s/3\rfloor -1$, which is impossible for $s\geq 3$.  

We finally need to verify the cases $s=3,4$. It is again sufficient to prove that if $H$ contains edges of weight 2 or 3 then $w(H)+|H_\Gamma|\leq t_3(s)+s-1$, thus we need to show that $w(H)+|H_\Gamma|$ is at most $5$ if $s=3$ and at most $8$ if $s=4$.

We note that argument in Case (i) above implies that if $H$ contains an edge of weight 3 then $|H_\Gamma|\leq s-2$ and $w(H)\leq 3+3\binom{s-2}{2}$, so if $s=3$ then $w(H)+|H_\Gamma|\leq 4$ and if $s=4 $ then  $w(H)+|H_\Gamma|\leq 8$ so the result holds. So we may suppose there are no edges of weight 3.

Now let $xy$ be an edge of weight $2$. Using the fact that $xy$ is not in any triangles and Lemma \ref{struct1:lem} (viii) and (ix) we find that for $s=3$ we have $w(H)+|H_\Gamma|\leq 2+3-|H_\Gamma|$, while for $s=4$ we have $w(H)+|H_\Gamma|\leq 2+ 6-|H_\Gamma|$, so the result holds. 
\qed

Finally we need to establish our two stuctural lemmas.

\emph{Proof of Lemma \ref{ax:lem}:} In each case we describe a labelling of the vertices of the given configuration to show that if it is present then $G$ is not $\F$-free.
\begin{itemize}
\item[($F_6$-1)] $a=1$, $b=5$, $c=6$, $x=2$, $y=3$, $z=4$.
\item[($F_6$-2)] $a=3$, $b=4$, $c=5$, $x=1$, $y=2$, $z=6.$
\item[($F_6$-3)] $a=1$, $b=2$, $c=3$, $x=4$, $y=5$, $z=6$.
\item[($F_6$-4)] $a=1$, $b=3$, $x=2$, $y=4$, $z=5$, $w=6$.
\item[($F_6$-5)] $a=5$, $b=1$, $c=3$, $x=4$, $y=2$, $z=6$.
\item[($K_4^-$-1)] $a=1$, $x=2$, $y=3$, $z=4$.
\item[($K_4^-$-2)] $a=3$, $b=4$, $x=1$, $y=2$.
\item[($K_4^-$-3)] $a=1$, $b=2$, $x=3$, $y=4$. \qed
\end{itemize}
\emph{Proof of Lemma \ref{struct1:lem}:}
We will make repeated use of Lemma \ref{ax:lem}.

(i) This follows immediately from $(F_6$-$1)$ and $(K_4^-$-$1)$.

(ii) This follows immediately from (i): if $uvwx$ is a copy of $K_4$ then we may suppose $l(uv)=a,l(uw)=b,l(vw)=c$, thus $l(ux)=c$ (otherwise (i) would be violated) continuing we see that $uvwx$ must be rainbow.

(iii) This follows immediately from (ii): if $xyzuv$ is a copy of $K_5$ then by (ii) we may suppose that $l(xy),l(xz),l(xu),l(xv)$ are all distinct single labels from $\{a,b,c\}$ but this is impossible since there are only 3 labels in total.

(iv) This follows immediately from  $(F_6$-$2)$ and $(K_4^-$-$2)$.

(v) If $x$ is in a $K_4$ then by (ii) it lies in edges with labels $a,b,c$, so $(F_6$-$3)$ implies that $x\not\in \Gamma_{abc}$.

(vi) If $x\in \Gamma_{abc}$, say $x\in \Gamma_{ab}$, and $y\in V_{abc}^4$ with $xy\in E(L_{abc})$ then $(F_6$-$3)$ implies that $l(xy)\neq c$, while $(F_6$-$4)$ implies that $l(xy)\neq a,b$ (since there are $t,u,v,w$ such that $l(yt)=b,l(tu)=a$ and $l(yv)=a,l(vw)=b$). 

(vii) This follows immediately from the fact that all $v\in V_{abc}^4$ meet edges with labels $a,b,c$ and $(F_6$-$2)$.

(viii) ($F_6$-5) implies that $\Gamma_{bc}=\emptyset$. If $xz\in E(L_{abc})$ then 
($F_6$-3) implies that $l(xz)=a$. Now $(K_4$-3) implies that $z\not\in \Gamma_{ac}$ while $(F_6$-3) implies that $z\not\in \Gamma_{bc}$. Hence $z\not\in \Gamma_{abc}$. Similarly if $yz\in E(L_{abc})$ then $l(yz)=b$ and $z\not\in \Gamma_{abc}$.

(ix) If $x\in\Gamma_{abc}$ or $y\in\Gamma_{abc}$ then this follows directly from (viii) so suppose that $x,y\not\in \Gamma_{abc}$, $l(xy)=ab$ and $|l(xz)|=2$. In this case, $(F_6$-2) implies that $l(xz)=ab$ so $(K_4$-2) implies that $z\in \Gamma_{ac}\cup\Gamma_{bc}$. But then $(F_6$-3) is violated. Hence $|l(xz)|\leq 1$. \qed

\section{Conclusion}
Many Tur\'an-type results have associated ``stability'' versions, and we were able to obtain such a result. For reasons of length we state it without proof.
\begin{thm}\label{thm:stabilityfsix}For any $\epsilon>0$ there
exist $\delta>0$ and $n_{0}$ such that the following holds: if
$H$ is an $\mathcal{F}$-free 3-graph of order $n\geq n_{0}$ with
at least $\left(1-\delta\right)s_3(n)$ edges, then there is a partition
of the vertex set of $H$ as $V(H)=U_{1}\cup U_{2}\cup U_{3}$ so
that all but at most $\epsilon n^{3}$ edges of $H$ have one vertex
in each $U_{i}$.\end{thm}


\begin{thebibliography}{}
\bibitem{BS} W. G. Brown and M. Simonovits, \emph{Digraph extremal problems, hypergraph extremal functions, and the densities of graph structures}, Disc. Math. {\bf 48} 147--162 (1984).
\bibitem{BT} R. Baber and J. Talbot, \emph{New Tur\'an densities for 3-graphs}, Electron. J. Combin. {\bf 19 (2)} P22 (2012). 
\bibitem{B} B. Bollob\'as, \emph{Three-graphs without two triples whose symmetric difference is contained in a third}, Disc. Math. {\bf 8} 21--24, (1974).
\bibitem{ESi} P. Erd\H os and M. Simonovits, \emph{A limit theorem in graph theory}, Studia Sci. Math. Hung. Acad. {\bf 1} 51--57, (1966).
\bibitem{ESt} P. Erd\H os and A.H. Stone, On the structure of linear graphs, Bull. Amer. Math. Soc. 52 (1946) 1087--1091.
\bibitem{FF} P. Frankl and Z. F\"uredi, \emph{A new generalization of the Erd\H os--Ko--Rado theorem}, Combinatorica {\bf 3} 341--349, (1983).
\bibitem{KNS} G. Katona, T. Nemetz and M. Simonovits, \emph{On a problem of Tur\'an in the theory of graphs}, Mat. Lapok {\bf 15}, 228--238, (1964).
\bibitem{KM} P. Keevash and D. Mubayi, \emph{Stability theorems for cancellative hypergraphs}, J. Combin. Theory Ser. B {\bf 92} 163--175 (2004).
\bibitem{M} V. W. Mantel, \emph{Problem 28}, Wiskundige Opgaven 10, (1907), 60--61. 
\bibitem{RF} A. A. Razborov, \emph{Flag Algebras}, Journal of Symbolic Logic, {\bf 72} (4) 1239--1282, (2007).
\bibitem{R4} A. A. Razborov, \emph{On $3$-hypergraphs with forbidden $4$-vertex configurations},in  SIAM J. Disc. Math. {\bf 24}, (3) 946--963 (2010).
\bibitem{T2} P. Tur\'an, On an extremal problem in graph theory, Mat. Fiz. Lapok 48 (1941) 436-452 [in Hungarian].
\end{thebibliography}
\end{document}